\documentclass[a4paper,12pt]{article}
\usepackage{amsmath,amsthm,amssymb}
\usepackage{texdraw}

\theoremstyle{plain}
\newtheorem{theorem}{Theorem}[section]
\newtheorem{lemma}[theorem]{Lemma}
\newtheorem{corollary}[theorem]{Corollary}
\theoremstyle{remark}
\newtheorem{example}[theorem]{Example}

\DeclareMathOperator{\centr}{center}
\DeclareMathOperator{\exc}{Exc}
\DeclareMathOperator{\lcm}{lcm}

\begin{document}
\title{Non-rational divisors over non-degenerate cDV-points}
\author{D.~A. Stepanov\thanks{The work was partially supported by 
Russian Foundation of Basic Research project no. 02-01-00441, and 
Grant of Leading Scientific Schools no. 489.2003.1.}}
\date{}

\maketitle

\begin{abstract}
Let $(X,o)$ be a 3-dimensional terminal singularity of type $cD$
or $cE$ defined in $\mathbb{C}^4$ by an equation non-degenerate with
respect to its Newton diagram. We show that there is not more than
1 non-rational divisor $E$ over $(X,o)$ with discrepancy $a(E,X)=1$.
We also describe all blowups $\sigma$ of $(X,o)$ such that
$E=\exc(\sigma)$ is non-rational and $a(E,X)=1$.
\end{abstract}

\section{Introduction}\label{S:intro}
In this paper, we study resolutions of 3-dimensional Gorenstein 
terminal singularities. The definition, the classification, and 
the basic properties of terminal singularities can be found in 
\cite{RY}. We need some concepts from toric geometry 
(see \cite{Danilov}); in particular, Varchenko-Hovanski\u{\i} embedded 
toric resolution is very important (see \cite{Varchenko}). In the 
sequel, all singularities are 3-dimensional and defined over the 
field $\mathbb{C}$ of complex numbers.

Let $\pi\colon Y\to X$ be a resolution of a terminal variety
(singularity) $X$, and let $E_i$ be a prime exceptional divisor of
$\pi$. It was proved in \cite{C3f}, 2.14 that $E_i$ is a
birationally ruled surface. Moreover, if $(X,o)$ is a germ of
singularity of type $cA/r$, $r\geq 1$, then all $E_i$ with
discrepancies $a(E_i,X)\leq 1$ and $\centr_X(E_i)=o$ are rational;
this was shown by Yu. G. Prokhorov in \cite{E1blowups}, 2.4. The 
meaning of the condition $a(E_i,X)<1$ is that exceptional divisors 
with this property appear on every resolution of $X$; exceptional 
divisors with discrepancy $a\leq 1$ and center at $o$ appear on every
divisorial resolution (here we identify two divisors over $X$ if they
give the same discrete valuation of the field $k(X)$).

However, there are non-rational divisors with $a\leq 1$ for other
types of terminal singularities. Some examples are given in section
\ref{S:examples}. In this paper, we describe non-rational exceptional
divisors with discrepancy $a=1$ and center at $o$, where $(X,o)$ is
a Gorenstein terminal singularity of type different from $cA$; i.e.,
it is a $cD$- or $cE$-singularity. Note that since $X$ is Gorenstein,
$1$ is the least possible value for $a(E,X)$. We also assume that
$(X,o)$ is defined by an equation non-degenerate with respect to its
Newton diagram. Our results are in Theorems \ref{T:uniqness} and 
\ref{T:description}, section \ref{S:preliminaries}.

If we don't assume non-degeneracy, the description of non-rational
divisors needs a more precise analysis. Such an analysis was carried
out for the $cD$-points in \cite{cD}. We include the $cD$-case here
for completeness.

I am grateful to Yu. G. Prokhorov and S. A. Kudryavtsev for useful
discussions and valuable advice.

\section{Preliminaries and results}\label{S:preliminaries}
Suppose that $X\subset\mathbb{C}^4$ is a hypersurface, $0\in X$. Then
$0\in X$ is a \emph{compound Du Val point} (a \emph{$cDV$-point}) if
a general hyperplane section $H$ of $X$ through $0$ is a surface with
Du Val singularity. We say that $0\in X$ \emph{is of type} $cA_n$, 
$n\geq 1$, ($cD_n$, $n\geq 4$, $cE_6$, $cE_7$, $cE_8$) if $(H,0)$ is
a Du Val singularity of type $A_n$ ($D_n$, $E_6$, $E_7$, $E_8$
respectively). In an analytic neighborhood of the point $0$ the
singularity $X$ of type $cA_n$ ($cD_n$, $cE_6$, $cE_7$, $cE_8$) can
be given by an equation of the form
$$f(x,y,z)+tg(x,y,z,t)=0\,,$$
where $g$ is some power series and $f=x^2+y^2+z^{n+1}$ ($x^2+y^2z+
z^{n-1}$, $x^2+y^3+z^4$, $x^2+y^3+yz^3$, $x^2+y^3+z^5$ respectively).
3-dimensional Gorenstein terminal singularities are exactly the
isolated $cDV$-points (\cite{Pagoda}).

Let $f(x,y,z,t)$ be a power series with Newton diagram $\Gamma(f)$.
Recall that $f$ is \emph{non-degenerate with respect to its Newton
diagram} (we shall simply say that $f$ is \emph{non-degenerate}) if
for every face $\rho$ of the diagram $\Gamma(f)$ the hypersurface
$$\{f_\rho(x,y,z,t)=0\}\subset(\mathbb{C}^*)^4$$
is smooth (see \cite{Varchenko}, \S2, \S5). Here 
$$f_\rho=\sum_{m\in\rho}a_mx^{m_1}y^{m_2}z^{m_3}t^{m_4}
\text{ if }
f=\sum_{\substack{m=\\ (m_1,m_2,m_3,m_4)}}
a_mx^{m_1}y^{m_2}z^{m_3}t^{m_4}\,.$$ 
\begin{lemma}\label{L:equations}
Suppose that $(X,0)\subset(\mathbb{C}^4,0)$ is a germ of an isolated
$cDV$-point of type different from $cA$. Then after a suitable 
coordinate change the singularity $(X,0)\subset(\mathbb{C}^4,0)$ is 
defined by one of the following equations.
\begin{equation}\label{E:cD}
x^2+y^2z+z^{n-1}+a_1t^{b_1}+a_2zt^{b_2}+\dots+a_{n-1}z^{n-2}t^{b_{n-1}}+
a_nyt^{b_n}+(\dots)=0\,,
\end{equation}
if $X$ is of type $cD_n$. Here $a_i\in\mathbb{C}$, $1\leq i\leq n$,
and for every $i$, $1\leq i\leq n-1$, we have $i-1+b_i\geq n-1$.
\begin{equation}\label{E:cE6}
x^2+y^3+z^4+a_1t^{b_1}+a_2zt^{b_2}+a_3z^2t^{b_3}+a_4yt^{b_4}+
a_5yzt^{b_5}+a_6yz^2t^{b_6}+(\dots)=0\,,
\end{equation}
if $X$ is of type $cE_6$. Here $i-1+b_i\geq 4$ for $i=1,2,3$.
\begin{align}\label{E:cE7}
x^2+y^3+& yz^3+a_1t^{b_1}+a_2zt^{b_2}+\dots+a_kz^{k-1}t^{b_k}+
a_{k+1}z^k+\\
& a_{k+2}yt^{b_{k+1}}+a_{k+3}yzt^{b_{k+2}}+(\dots)=0\,,\notag
\end{align}
if $X$ is of type $cE_7$. Here $k$ is some integer $\geq 5$.
\begin{align}\label{E:cE8}
x^2+y^3+z^5+a_1t^{b_1}+a_2zt^{b_2}+a_3z^2t^{b_3}+a_4z^3t^{b_4}+\\
a_5yt^{b_5}+a_6yzt^{b_6}+a_7yz^2t^{b_7}+a_8yz^3t^{b_8}+(\dots)=0\,,
\notag
\end{align}
if $X$ is of type $cE_8$. Here $i-1+b_i\geq 5$ for $i=1,2,3,4$.

In every equation, the part $(\dots)$ does not affect its Newton
diagram.
\end{lemma}
\begin{proof}
We prove the lemma for $cE_6$-singularities. Other cases can be done 
in a similar way.

So, let $X$ be given by the equation
$$f=x^2+y^3+z^4+tg(x,y,z,t)=0\,.$$
Put the monomials of $tg$ in the increasing order with respect to 
grlex-ordering (see \cite{Cox}, 2 \S2). Suppose that $a_{\alpha\beta\gamma
\delta}x^\alpha y^\beta z^\gamma t^\delta$ is the first monomial in
$tg$ containing the variable $x$. If $\alpha=1$, then we make the
coordinate change $x\leftarrow x-\frac{a_{\alpha\beta\gamma\delta}}
{2}y^\beta z^\gamma t^\delta$, 
$y\leftarrow y$, $z\leftarrow z$, $t\leftarrow t$.
If $\alpha\geq 2$, then we substitute
$x\leftarrow x/\sqrt{1+a_{\alpha\beta\gamma\delta}x^{\alpha-2}
y^\beta z^\gamma t^\delta}$, $y\leftarrow y$, etc.
In both cases we remove the monomial
$a_{\alpha\beta\gamma\delta}x^\alpha y^\beta z^\gamma t^\delta$
but the previous part of $tg$ is unchanged. Repeating this procedure
we increase the number 
$$\min\{\deg x^\alpha y^\beta z^\gamma 
t^\delta\,|\,x^\alpha y^\beta z^\gamma t^\delta\in tg,
 \alpha\geq 1\}\,.$$

Recall that the given singularity $(X,0)$ is isolated. Therefore it 
is finitely determined (\cite{Hir}, Theorem 3.3)
and hence we can reduce the equation of $X$ to the form
$$x^2+y^3+z^4+tg(y,z,t)=0\,.$$
(We should write $g'$ instead of $g$ but we hope that there will be
no confusion.) 

Similarly, using the terms $y^3$ and $z^4$ we can remove monomials 
with $y^\beta$, $\beta\geq 2$, and $z^\gamma$, $\gamma\geq 3$ from $tg$
and reduce the equation to the form \eqref{E:cE6}. If there is a
monomial $z^\gamma t^\delta$ in $f$, $\gamma+\delta\leq 3$, then $f$ 
defines a singularity of type $cD_4$ or $cA$; this proves the 
condition $i-1+b_i\geq 4$ for $i=1,2,3$. 
The fact that terms in $(\dots)$ do not affect the Newton diagram 
follows from the construction.
\end{proof}

\begin{theorem}\label{T:uniqness}
Let $(X,0)$ be a terminal point of type $cD$ or $cE$ defined in
$\mathbb{C}^4$ by one of equations \eqref{E:cD}--\eqref{E:cE8}; also,
assume that this equation is non-degenerate with respect to its
Newton diagram. Then for every resolution $\pi\colon Y\to X$ there
exists at most one non-rational exceptional divisor $E$ of $\pi$
such that $a(E,X)=1$ and $\centr_X(E)=0$.
\end{theorem}
\begin{theorem}\label{T:description}
In the conditions of Theorem \ref{T:uniqness} assume that $E$ is the 
non-rational divisor over $X$ such that $a(E,X)=1$ and 
$\centr_X(E)=0$. Then $E$ is (birational to) the exceptional divisor
of the weighted blowup $\sigma_w$, where the weight $w$ is described 
below.
\begin{description}
\item[(i)] If $X$ is of type $cD_n$, then \\
1) if $n=2k$, $k\geq 2$, then $w=(k,k-1,1,1)$; \\
2) if $n=2k+1$, $k\geq 2$, then $w=(k,k,1,1)$. \\
In both cases $E$ is birational to the surface $C\times\mathbb{P}^1$,
where $C$ is a hyperelliptic curve of genus $g\leq k-1$.
\item[(ii)] If $X$ is of type $cE_6$, then one of the following
holds: \\
1) $w=(2,2,1,1)$; \\
2) $w=(3,2,2,1)$; \\
3) $w=(4,3,2,1)$. \\
In all cases $E\simeq C\times\mathbb{P}^1$, where $C$ is a curve of 
genus $1$.
\item[(iii)] If $X$ is of type $cE_7$, then one of the following
holds: \\
1) $w=(3,2,1,1)$; 
2) $w=(4,3,2,1)$; \\
3) $w=(5,3,2,1)$; 
4) $w=(6,4,3,1)$. \\
In cases 1), 2), 4) $E\simeq C\times\mathbb{P}^1$, where $C$ is a 
curve of genus 1. In case 3) $g(C)\leq 3$ and $C$ can be 
non-hyperelliptic.
\item[(iv)] If $X$ is of type $cE_8$, then one of the following 
holds: \\
1) $w=(3,2,2,1)$; 
2) $w=(4,3,2,1)$; 
3) $w=(5,3,2,1)$; \\
4) $w=(6,4,3,1)$; 
5) $w=(7,5,3,1)$; 
6) $w=(8,5,3,1)$; \\
7) $w=(9,6,4,1)$; 
8) $w=(12,8,5,1)$. \\
In case 6) $E\simeq C\times\mathbb{P}^1$, where $g(C)\leq 4$ and
$C$ can be non-hyperelliptic. In other cases $g(C)=1$.
\end{description}
\end{theorem}
For the \emph{proof}, see section \ref{S:proof}.

Let a hypersurface $X\subset \mathbb{C}^n$ be defined by a
non-degenerate function $f(x_1,\dots,x_n)$. Suppose that $0\in X$
is an isolated singularity. Then there exists an embedded toric
resolution for $X\subset\mathbb{C}^n$, i.e., a toric morphism
$\pi\colon\widetilde{\mathbb{C}^n}\to\mathbb{C}^n$ such that
$\widetilde{\mathbb{C}^n}$ and the proper transform $\widetilde X$
of $X$ are nonsingular, and the exceptional locus of 
$\pi|_{\widetilde X}$ is a divisor with normal crossings 
(\cite{Varchenko}, \S9, \S10). The morphism $\pi$ is determined by a
certain subdivision $\Sigma$ of the positive octant 
$\mathbb{R}^{n}_{\geq 0}$. Exceptional divisors $E_i$ of $\pi$
correspond to 1-dimensional cones $\tau_i\subset\mathbb{R}^{n}_{>0}$
of $\Sigma$. Fix some prime exceptional divisor $E_\tau$ and the
corresponding cone $\tau$, take a primitive vector $w=
(w_1,\dots,w_n)$ along $\tau$ ($w\in\mathbb{Z}^n\subset
\mathbb{R}^n$), and let $E_\tau|_{\widetilde X}=\sum_j m_jE_j$. Also
denote by $\Gamma(f)$ the Newton diagram of $f$; $\Gamma(f)\subset
(\mathbb{R}^n)^*$, where the space $(\mathbb{R}^n)^*$ is dual to
$\mathbb{R}^n$; we denote the corresponding pairing by 
$\langle \cdot,\cdot\rangle$. 

Now we want to calculate the discrepancy $a(E_j,X)$.
\begin{lemma}\label{L:discrepancy}
$a(E_j,X)=m_j(w_1+w_2+\dots+w_n-1-w(f))$, where $w(f)=
\min\{\langle w,v\rangle\,|\,v\in\Gamma(f)\}$.
\end{lemma}
\begin{proof}
According to \cite{Varchenko}, \S10, there is an affine neighborhood
$U\simeq\mathbb{C}^n$ of $E_\tau$ in $\widetilde{\mathbb{C}^n}$. If
$y_1,\dots,y_n$ are coordinates in $U$ such that the equation $y_1=0$
defines $E_\tau\cap U$, then $\pi|_U\colon U\to\mathbb{C}^n$ is given
by
$$x_1=y_{1}^{w_1}y_{2}^{a^{2}_{1}}\dots y_{n}^{a^{n}_{1}}\,,$$
$$\dots\dots$$
$$x_n=y_{1}^{w_n}y_{2}^{a^{2}_{n}}\dots y_{n}^{a^{n}_{n}}$$
for some $a^i=(a^{i}_{1},\dots,a^{i}_{n})\in\mathbb{Z}^{n}_{\geq 0}$.
Now to prove the lemma we need only to lift the differential form
$dx_1\wedge\dots\wedge dx_n$ to $U$ and apply the adjunction formula.
\end{proof}
\begin{corollary}\label{C:discrepancy}
If $X$ is terminal Gorenstein and $a(E_j,X)=1$, then
$$w_1+w_2+\dots+w_n-1-w(f)=1\,.$$
\end{corollary}

Now suppose that $0\in X\subset\mathbb{C}^4$ is a terminal 
$cDV$-point and let $\sigma_w\colon\widetilde{\mathbb{C}^4}\to
\mathbb{C}^4$ be the weighted blowup with weight $w=
(w_1,w_2,w_3,w_4)$. We use the same notation $\sigma_w\colon
\widetilde X\to X$ for the restriction of $\sigma_w$ to the proper
transform $\widetilde X\subset\widetilde{\mathbb{C}^4}$ of $X$. In
this situation the exceptional divisor $E\subset\widetilde X$ of
$\sigma_w$ is a surface in the weighted projective space 
$\mathbb{P}(w_1,w_2,w_3,w_4)$.
\begin{lemma}\label{L:rationality}
Notation as above. Assume that $E\subset
\mathbb{P}(w_1,w_2,w_3,w_4)$ is irreducible
and has rational singularities. Then $E$ is rational.
\end{lemma}
\begin{proof}
Let $\pi\colon\widetilde E\to E$ be a resolution of $E$. By 
\cite{C3f}, 2.14, $E$ is birationally ruled. Thus $P_2(\widetilde E)=
h^0(2K_{\widetilde E})=0$. On the other hand, $E$ is a hypersurface
in $\mathbb{P}(w_1,w_2,w_3,w_4)$, hence $h^1(\mathcal{O}_E)=0$. Since
$E$ has rational singularities, we have 
$h^1(\mathcal{O}_{\widetilde E})=h^1(\mathcal{O}_E)=0$. Therefore
$\widetilde E$ is rational by Castelnuovo's criterion.
\end{proof}
For example, Lemma \ref{L:rationality} can be applied when $E$ is
quasismooth or when $\sigma_w$ is a plt-blowup (see \cite{blow-ups},
Definition 3).

\section{Proof of theorems \ref{T:uniqness} and \ref{T:description}}
\label{S:proof}
Take an embedded toric resolution $\pi\colon\widetilde X\to X$ of the
given non-degenerate terminal $cDV$-point $0\in X\subset\mathbb{C}^4$,
$X=\{f=0\}$. Let $\Sigma$ be the corresponding subdivision of
$\mathbb{R}^{4}_{\geq 0}$ and take a primitive vector 
$w\in\mathbb{Z}^{4}_{>0}$ along a 1-dimensional cone $\tau$ of 
$\Sigma$. As in Lemma \ref{L:discrepancy}, suppose $E_\tau\cap
\widetilde X=\cup_jE_j$. Since for any two subdivisions of
$\mathbb{R}^{4}_{\geq 0}$ there is a common subdivision, we see that
the divisors $E_j$ are in 1-to-1 correspondence with the exceptional
divisors of the weighted blowup $\sigma_w$. If $E'_j\subset
\exc(\sigma_w)$ corresponds to $E_j$, then $E'_j$ and $E_j$ are
birational.

The divisor $\exc{\sigma_w}$ is given in 
$\mathbb{P}(w_1,w_2,w_3,w_4)$ by the equation 
$$f_{\rho(w)}(x,y,z,t)=0\,,$$ 
where $f_{\rho(w)}$ is the part of $f$
corresponding to the face 
$$\rho(w)=\{v\in\Gamma(f)\,|\,\langle w,v\rangle=w(f)\}\,.$$

Now let us consider cases $cD_n$, $cE_6$, $cE_7$, $cE_8$ one by one.

\subsection{$cD_n$}\label{S:cD}
\begin{figure}[ht]
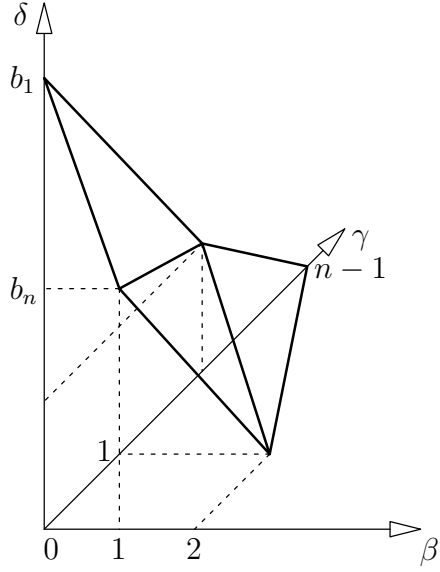

\centertexdraw{
\drawdim cm \linewd 0.02 
\htext(0 -0.4){$0$} 
\move(0 0) \avec(5 0) \htext(5 -0.5){$\beta$} 
\move(0 0) \avec(4 4) \htext(4.1 3.7){$\gamma$}
\move(0 0) \avec(0 7) \htext(-0.4 6.7){$\delta$}
\linewd 0.03
\move(3 1) \lvec(3.5 3.5) \lvec(2.1 3.8)
\lvec(0 6) \lvec(1 3.2) \lvec(2.1 3.8)
\lvec(3 1) \lvec(1 3.2)
\linewd 0.02 \lpatt(0.05 0.1)
\lvec(1 0) \move(1 3.2) \lvec(0 3.2)
\move(3 1) \lvec(2 0) \move(3 1) \lvec(1 1)
\move(2.1 3.8) \lvec(2.1 2.1) \move(2.1 3.8) \lvec(0 1.7)
\htext(0.9 -0.4){$1$} \htext(1.9 -0.4){$2$} 
\htext(3.6 3.3){$n-1$} \htext(0.7 0.9){$1$}
\htext(-0.45 5.8){$b_1$} \htext(-0.45 3){$b_n$}
}
\caption{Newton diagram for $cD$-singularity}\label{Fig:cD}
\end{figure}

Suppose that $0\in X$ is of type $cD_n$, $n\geq 4$. Then it is given
in $\mathbb{C}^4$ by equation \eqref{E:cD}:
$$f=x^2+y^2z+z^{n-1}+a_1t^{b_1}+a_2zt^{b_2}+\dots+a_{n-1}z^{n-2}t^{b_{n-1}}+
a_nyt^{b_n}+(\dots)=0\,.$$
By our assumption, $f$ is non-degenerate.
An example of a Newton diagram for $f(0,y,z,t)$ is presnted in 
Figure \ref{Fig:cD}. Note that \eqref{E:cD} contains only one 
monomial with the variable $x$ ($x^2$), thus $\Gamma(f)$ is 
completely determined by the Newton diagram for $f(0,y,z,t)$.

Assume that $E$ is a non-rational divisor over $X$ such that
$a(E,X)=1$ and $\centr_X(E)=0$. Then we can consider $E$ as one of
exceptional divisors of a certain weighted blowup $\sigma_w$. The
exceptional divisor of $\sigma_w$ corresponds to some face $\rho$
of $\Gamma(f)$. It is clear that 0-dimensional faces produce only
rational divisors. If $\rho$ is a 1-dimensional face, then 
$f_\rho$ is equal either to a quasihomogeneous polynomial of
$z$ and $t$, or to a polynomial containing only two monomials. In
both cases, $f_\rho$ defines rational surfaces. Thus $\dim\rho=2$
or $3$.

If $\rho$ does not contain any of the monomials $x^2$,
$y^2z$, $z^{n-1}$, then $f_\rho=a_nyt^{b_n}+h(z,t)$, $a_n\ne 0$.
Again $\{f_\rho=0\}$ is rational. Thus we may assume that $x^2\in
\rho$ or $y^2z\in\rho$.

Now let us use the condition $a(E,X)=1$. Let $L$ be the hyperplane
in $(\mathbb{R}^4)^*$ such that $L\supset\rho$ and $w$ is normal
to $L$. Write the equation of $L$ in the form
$$\frac{\alpha}{a}+\frac{\beta}{b}+\frac{\gamma}{c}+\frac{\delta}{d}
=1\,,$$
where $\alpha$, $\beta$, $\gamma$, $\delta$ are the coordinates in
$(\mathbb{R}^4)^*$. If $m=\lcm(a,b,c,d)$, then $m=w(f)$ and 
$m/a=w_1$, $m/b=w_2$, $m/c=w_3$, $m/d=w_4$. In our case, \\
(i) ($a=2$, $\frac{2}{b}+\frac{1}{c}\geq 1$, $c\leq n-1$) or 
($a<2$, $\frac{2}{b}+\frac{1}{c}=1$, $c\leq n-1$); \\
(ii) $\frac{m}{a}+\frac{m}{b}+\frac{m}{c}+\frac{m}{d}-1-m=1$ (see
Corollary \ref{C:discrepancy}).
\begin{lemma}\label{L:cDnum}
Let $a$, $b$, $c$, $d\in\mathbb{Q}_{>0}$ satisfy conditions (i) and
(ii). Then $(a,b,c,d)$ is one of the following: \\
1)$(\frac{2k-1}{k},\frac{2k-1}{k-1},2k-1,2k-1)$; here $n=2k$, 
$k\geq 2$; \\
2)$(2,2,2k,2k)$; here $n=2k+1$, $k\geq 2$; \\
3)$(2,\frac{2k-1}{k-1},k,2k)$, $k\leq n-1$.
\end{lemma}
\begin{proof}
Case 1: $a=2$. In this case $m=2k$, $k\geq 2$.

We have $2k(1/2+1/b+1/c+1/d)-1-2k=1$. Using $2/b+1/c\geq 1$, we get
$$\frac{2k}{2c}+\frac{2k}{d}\leq 2\,.$$
But $2k/d\in\mathbb{Z}_{>0}$, hence $2k/d=1$, 
$2k/c\leq 1$. Also $2k/c\in\mathbb{Z}$; thus 
$$\frac{2k}{2c}=\frac{1}{2}\text{ or }\frac{2k}{2c}=1\,.$$

Let $\frac{2k}{2c}=\frac{1}{2}$, $c=2k=m=d$. Then $c=n-1$, $b=2$. 
We obtain case 2).

Now suppose $\frac{2k}{2c}=1$, $k=c$. Here we obtain 
$b=\frac{2k}{k-1}$ and
$$(a,b,c,d)=(2,\frac{2k}{k-1},k,2k)\,,$$ 
where $k\leq n-1$. This is case 3).

Case 2: $a<2$, $2/b+1/c=1$. In this case, from (ii) we get 
$$\frac{m}{2c}+\frac{m}{d}<2\,.$$
But $\frac{m}{d}\in\mathbb{Z}_{>0}$, 
$\frac{m}{2c}\in\frac{1}{2}\mathbb{Z}_{>0}$, therefore $m=d$, 
$\frac{m}{2c}=\frac{1}{2}$, hence $m=c=d$. Thus $m=n-1$, $b=
\frac{2(n-1)}{n-2}$. Substituting $\frac{2(n-1)}{n-2}$, $n-1$, $n-1$
for $b$, $c$, $d$ in (ii), we obtain 
$$\frac{n-1}{a}+\frac{n-1}{2}=\frac{1}{2}\,.$$
Since $\frac{n-1}{a}\in\mathbb{Z}$, we see that $n-1$ is odd, i.e., 
$n-1=2k-1$. We get case 1).
\end{proof}

Cases 1) and 2) of Lemma \ref{L:cDnum} correspond to blowups (i),
1) and 2) of Theorem \ref{T:description}. Case 3) gives us the
blowup $\sigma=(k,k-1,2,1)$. But $\sigma$ satisfies the conditions
of \cite{IP}, Proposition 3.2, hence it is a plt-blowup of $X$ and
it follows from Lemma \ref{L:rationality} that its exceptional 
divisor is rational. In other cases the exceptional divisor can be 
non-rational (see examples in section \ref{S:examples}), but in case
1) $n$ is even and in case 2) $n$ is odd, therefore the non-rational
divisor is unique.

For the weighted blowup $\sigma=(k,k-1,1,1)$ (Theorem 
\ref{T:description}, (i) 1)), we get that the exceptional divisor $E$
is defined in $\mathbb{P}(k,k-1,1,1)$ by the equation
$$y^2z+z^{2k-1}+h(y,z,t)=0\,,$$
where $h=yh_1(z,t)+h_2(z,t)$ is a quasihomogeneous polynomial of 
degree $2k-1$ with respect to the weights $(k-1,1,1)$.  We see that 
$E$ is a cone over a hyperelliptic curve $C$ of degree $2k-1$ in
$\mathbb{P}(k-1,1,1)$. It follows easily that genus $g(C)\leq k-1$.

For the weighted blowup $\sigma=(k,k,1,1)$ (Theorem
\ref{T:description}, (i) 2)), the exceptional divisor $E$ is given
in $\mathbb{P}(k,k,1,1)$ by the equation 
$$x^2+h(z,t)=0\,,$$
where $h$ is a homogeneous polynomial of  degree $2k$. Again $E$ is
a cone over a hyperelliptic curve of genus $g\leq k-1$. This proves
Theorems \ref{T:uniqness} and \ref{T:description} for
$cD$-singularities.

\subsection{$cE_6$}\label{S:cE6}
\begin{figure}[ht]
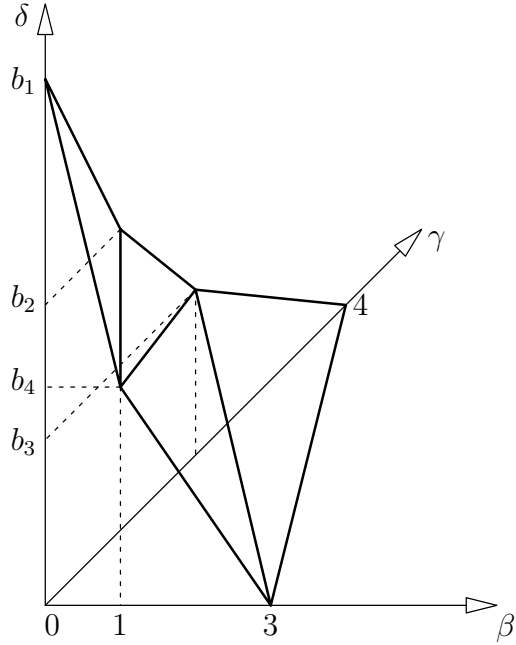

\centertexdraw{
\drawdim cm \linewd 0.02 
\htext(0 -0.4){$0$} 
\move(0 0) \avec(6 0) \htext(6 -0.5){$\beta$} 
\move(0 0) \avec(5 5) \htext(5.1 4.7){$\gamma$}
\move(0 0) \avec(0 8) \htext(-0.4 7.7){$\delta$}
\linewd 0.03
\move(3 0) \lvec(4 4) 
\lvec(2 4.2) \lvec(1 5) \lvec(0 7) 
\lvec(1 2.9) \lvec(3 0) 
\lvec(2 4.2) \lvec(1 2.9) \lvec(1 5)
\linewd 0.02 \lpatt(0.05 0.1)
\move(2 4.2) \lvec(2 2) \move(2 4.2) \lvec(0 2.2)
\move(1 5) \lvec(0 4) 
\move(1 2.9) \lvec(1 0) \move(1 2.9) \lvec(0 2.9)
\htext(0.9 -0.4){$1$} \htext(2.9 -0.4){$3$} 
\htext(4.1 3.85){$4$} 
\htext(-0.45 6.8){$b_1$} \htext(-0.45 3.9){$b_2$} 
\htext(-0.45 2){$b_3$} \htext(-0.45 2.8){$b_4$}
}
\caption{Newton diagram for $cE_6$-singularity}\label{Fig:cE6}
\end{figure}

Suppose that $0\in X$ is of type $cE_6$. Then it is defined in
$\mathbb{C}^4$ by equation \eqref{E:cE6}
$$f=x^2+y^3+z^4+a_1t^{b_1}+a_2zt^{b_2}+a_3z^2t^{b_3}+$$
$$a_4yt^{b_4}+a_5yzt^{b_5}+a_6yz^2t^{b_6}+(\dots)=0\,.$$
By our assumption, $f$ is non-degenerate. 
An example of a Newton diagram for $f(0,y,z,t)$ is presented
in Figure \ref{Fig:cE6}.
As above, $\Gamma(f)$ is completely determined by the Newton diagram
for $f(0,y,z,t)$.

Let $E$ be a non-rational divisor over $X$ such that $a(E,X)=1$ and
$\centr_X(E)=0$. Assume that $E$ is an exceptional divisor of a
certain weighted blowup $\sigma_w$; let $\rho$ be the corresponding
face of $\Gamma(f)$. The reader will easily prove that 0- and
1-dimensional faces produce rational divisors, so that we must only 
consider the cases $\dim\rho=2$ and $\dim\rho=3$.

If $\rho$ does not contain any of the monomials $x^2$, $y^3$,
$z^4$, then it must contain at least one of the monomials $yt^{b_4}$,
$yzt^{b_5}$, $yz^2t^{b_6}$ and thus $E$ is rational.

Let $L\subset(\mathbb{R}^4)^*$ be a hyperplane such that $L\supset
\rho$ and $w$ is normal to $L$. Write the equation of $L$ in the
form
$$\frac{\alpha}{a}+\frac{\beta}{b}+\frac{\gamma}{c}+\frac{\delta}{d}
=1\,,$$
where $\alpha$, $\beta$, $\gamma$, $\delta$ are the coordinates in
$(\mathbb{R}^4)^*$. If $m=\lcm(a,b,c,d)$, then $m=w(f)$ and $w_1=
m/a$, $w_2=m/b$, $w_3=m/c$, $w_4=m/d$. In our case, \\
(i) ($a=2$, $b\leq 3$, $c\leq 4$) or ($a<2$, $b=3$, $c\leq 4$)
or ($a<2$, $b<3$, $c=4$). \\
From the condition $a(E,X)=1$, by Corollary \ref{C:discrepancy} we
obtain \\
(ii) $m(\frac{1}{a}+\frac{1}{b}+\frac{1}{c}+\frac{1}{d})-1-m=1$. 
\begin{lemma}\label{L:cE6num}
Let $a$, $b$, $c$, $d\in\mathbb{Q}_{>0}$ satisfy conditions (i) and
(ii). Then $(a,b,c,d)$ is one of the following: \\
1) $(2,2,4,4)$; \\
2) $(2,3,3,6)$; \\
3) $(2,8/3,4,8)$; \\
4) $(2,3,4,12)$.
\end{lemma}
\begin{proof}
Case 1: $a=2$. Here $m=2k$. We have
$$\frac{1}{a}+\frac{1}{b}+\frac{1}{c}\geq \frac{1}{2}+\frac{1}{3}+
\frac{1}{4}=\frac{13}{12},.$$
Thus $m(1/a+1/b+1/c+1/d)-1-m\geq m(\frac{1}{12}+\frac{1}{d})-1$. 
From (ii) it follows that
$$\frac{m}{12}+\frac{m}{d}\leq 2\,.$$
Since $m/d\in\mathbb{Z}_{>0}$, we obtain $m/d=1$, $m\leq 12$. On the
other hand, it is clear that $m\geq 4$.

Let $m=d=4$. Then from (ii) we get
$$4(\frac{1}{2}+\frac{1}{b}+\frac{1}{c}+\frac{1}{4}-1)=2\,,$$
$$\frac{4}{b}+\frac{4}{c}=3\,.$$
There is only one possibility $b=2$, $c=4$. This gives case 1).

Let $m=d=6$. Here we get $6/b+6/c=4$. The only possibility is
$b=3$, $c=3$. This gives us case 2).

Let $m=d=8$. Here we find $b=8/3$, $c=4$. This is case 3).

Let $m=d=10$. We get $10/b+10/c=6$. Since $b\leq 3$, $10/b\geq 4$;
since $c\leq 4$, $10/c\geq 3$, a contradiction.

Let $m=d=12$. The only possibility is $b=3$, $c=4$, i.e., case 4).

Case 2: $a<2$. In this case
$$m(\frac{1}{a}+\frac{1}{b}+\frac{1}{c}+\frac{1}{d})-1-m>
m(\frac{1}{12}+\frac{1}{d})-1\,.$$
Thus $\frac{m}{12}+\frac{m}{d}<2$, $m<12$, $d=m$.

Case 2.1: $b=3$. Here $m=3k$, i.e., $m=6$ or $9$.

Let $m=d=6$. We get
$$6(\frac{1}{a}+\frac{1}{3}+\frac{1}{c}+\frac{1}{6}-1)=2\,,$$
$$\frac{6}{a}+\frac{6}{c}=5\,.$$
Since $6/a\in\mathbb{Z}$ and $a<2$, we have $6/a=4$, hence $c=6$.
This contradicts (i).

Let $m=d=9$. Here $9/a+9/c=7$. In the same way we obtain a 
contradiction.

Case 2.2: $b<3$, $c=4$. In this case $m=4k$, i.e., $m=4$ or $8$. In
both cases we easily come to a contradiction.
\end{proof}

Cases 1), 2), 3) of Lemma \ref{L:cE6num} correspond to blowups (ii)
1), 2), 3) of Theorem \ref{T:description} respectively. Case 4) gives
the blowup $\sigma=(6,4,3,1)$. It is a plt-blowup (\cite{IP}, 3.2) 
and its exceptional divisor is rational (Lemma \ref{L:rationality}).

Now we show that if a non-rational divisor $E$ with $a(E,X)=1$ exists,
then it is unique. Indeed, suppose that $E$ is the exceptional
divisor of the blowup $(2,2,1,1)$. Then $E$ is defined in
$\mathbb{P}(2,2,1,1)$ by the equation
$$x^2+z^4+a_1t^4+a_2zt^3+a_3z^2t^2=0\,,$$
where $a_1\ne 0$ or $a_2\ne 0$ (since $E$ is non-rational). It
follows that in this case exceptional divisors of the blowups
$(3,2,2,1)$ and $(4,3,2,1)$ are given by the equations $t^4=0$ or
$zt^3=0$ and hence they are rational.

Now suppose that the blowup $\sigma=(3,2,2,1)$ has a non-rational
exceptional divisor $E$. This means that the Newton diagram 
$\Gamma(f)$ lies above the hyperplane
$$\frac{\alpha}{2}+\frac{\beta}{3}+\frac{\gamma}{4}+
\frac{\delta}{6}=1$$
and $E$ is given in $\mathbb{P}(3,2,2,1)$ by the equation
$$x^2+y^3+a_1t^6+a_2zt^4+a_3z^2t^2+a_4yt^4+a_5yzt^2=0\,,$$
but since $E$ is non-rational, we have $a_2=a_3=a_5=0$ (otherwise
$\sigma$ is a plt-blowup). Here we see that exceptional
divisor of the blowup $(4,3,2,1)$ is rational.

In all cases, the non-rational divisor $E$ is a cone over a curve
of genus $1$. This proves Theorems \ref{T:uniqness} and
\ref{T:description} for $cE_6$-singularities.

\subsection{$cE_7$}\label{S:cE7}
Suppose that $0\in X$ is of type $cE_7$. Then it is defined in
$\mathbb{C}^4$ by equation \eqref{E:cE7}
$$f=x^2+y^3+yz^3+a_1t^{b_1}+a_2zt^{b_2}+\dots+a_kz^{k-1}t^{b_k}+
a_{k+1}z^k+$$
$$a_{k+2}yt^{b_{k+1}}+a_{k+3}yzt^{b_{k+2}}+(\dots)=0\,.$$
We assume that $f$ ia non-degenerate. 
As above, $\Gamma(f)$ is completely determined by the Newton diagram 
for $f(0,y,z,t)$.

Assume that a weighted blowup $\sigma_w$ of $X$ has a non-rational
exceptional divisor $E$ with discrepancy $1$ and let $\rho$ be the
corresponding face of $\Gamma(f)$. We can only consider the cases
$\dim\rho=2$ and $\dim\rho=3$.

If $\rho$ does not contain any of the monomials $x^2$, $y^3$, then
it must contain  the monomial $yt^{b_{k+2}}$, and thus $E$ is 
rational.

Using notation of the previous sections, we obtain 4 rational numbers
$a$, $b$, $c$, $d$ with the following properties. \\
(i) ($a=2$, $b\leq 3$, $1/b+3/c\geq 1$) or ($a<2$, $b=3$, 
$c\leq 9/2$); \\
(ii) $m(\frac{1}{a}+\frac{1}{b}+\frac{1}{c}+\frac{1}{d})-1-m=1$
($m=\lcm(a,b,c,d)$). 
\begin{lemma}\label{L:cE7num}
Let $a$, $b$, $c$, $d\in\mathbb{Q}_{>0}$ satisfy conditions (i) and
(ii). Then $(a,b,c,d)$ is one of the following: \\
1) $(2,2,6,6)$; 
2) $(2,3,3,6)$; 
3) $(2,8/3,4,8)$; \\
4) $(9/5,3,9/2,9)$; 
5) $(2,5/2,5,10)$; 
6) $(2,3,4,12)$; \\
7) $(2,14/5,14/3,14)$; 
8) $(2,3,9/2,18)$.
\end{lemma}
\begin{proof}
Case 1: $a=2$. Here $m=2k$.

Using the condition $1/b+3/c\geq 1$, from (ii) we get
$$\frac{4k}{3b}+\frac{2k}{3}+\frac{2k}{d}-k\leq 2\,.$$
Since $b\leq 3$, we have
$$\frac{4k}{9}+\frac{2k}{3}+\frac{2k}{d}-k\leq 2\,,$$
$$\frac{k}{9}+\frac{2k}{d}\leq 2\,.$$
But $2k/d\in\mathbb{Z}$, so that $d=2k$, $2k\leq 18$. The rest of
the proof is similar to the proof of Lemma \ref{L:cE6num}.
\end{proof}

Cases 1), 3), 4), 6) of Lemma \ref{L:cE7num} correspond to blowups
(iii) 1), 2), 3), 4) of Theorem \ref{T:description} respectively.
Other cases give us blowups with rational exceptional divisors. By
the same argument as in case $cE_6$, we can prove that if $E$ is
non-rational and given by one of the blowups (iii) 1), 2), 3), 4)
(Theorem \ref{T:description}), then other blowups have only rational 
exceptional divisors. Thus the non-rational divisor $E$ is unique.

In all cases, the non-rational divisor is a cone over a curve $C$.
It is not difficult to verify by a direct calculation that genus 
$g(C)$ satisfies the conditions of Theorem \ref{T:description}, (iii).

\subsection{$cE_8$}\label{S:cE8}
Suppose that $0\in X$ is of type $cE_8$. Then it is defined in 
$\mathbb{C}^4$ by equation \eqref{E:cE8}
\begin{align*}
f=x^2+y^3+z^5+a_1t^{b_1}+a_2zt^{b_2}+a_3z^2t^{b_3}+a_4z^3t^{b_4}+\\
a_5yt^{b_5}+a_6yzt^{b_6}+a_7yz^2t^{b_7}+a_8yz^3t^{b_8}+(\dots)=0\,.
\end{align*}
We assume that $f$ is non-degenerate. As above, $\Gamma(f)$ is
completely determined by the Newton diagram of $f(0,y,z,t)$. 
Using notation and argument of section \ref{S:cE6}, we obtain 4
rational numbers $a$, $b$, $c$, $d$ such that \\
(i) ($a=2$, $b\leq 3$, $c\leq 5$) or ($a<2$, $b=3$, $c\leq 5$) or
($a<2$, $b<3$, $c=5$); \\
(ii) $m(\frac{1}{a}+\frac{1}{b}+\frac{1}{c}+\frac{1}{d})-1-m=1$
($m=\lcm(a,b,c,d)$).
\begin{lemma}\label{L:cE8num}
Let $a$, $b$, $c$, $d\in\mathbb{Q}_{>0}$ satisfy conditions (i) and
(ii). Then $(a,b,c,d)$ is one of the following: \\
1) $(2,3,3,6)$; 
2) $(2,8/3,4,8)$; 
3) $(9/5,3,9/2,9)$; \\
4) $(2,3,4,12)$; 
5) $(2,14/5,14/3,14)$; 
6) $(15/8,3,5,15)$; \\
7) $(2,3,9/2,18)$; 
8) $(2,3,24/5,24)$;
9) $(2,3,5,30)$.
\end{lemma}
The \emph{proof} is similar to the proof of Lemma \ref{L:cE6num}.

Cases 1)--8) of Lemma \ref{L:cE8num} correspond to blowups (iv) 1)--8)
of Theorem \ref{T:description} respectively. The 9th case gives the
blowup $(15,10,6,1)$. This is a plt-blowup and its exceptional
divisor is rational.

By the same argument as in the previous sections, we can prove that 
if the nonrational divisor $E$ exists, then it is unique. In all cases 
$E$ is a cone over a curve in a weighted projective plane and we can
easily verify the rest of assertions of Theorem \ref{T:description}.

\section{Examples}\label{S:examples}
\begin{example}
Let $X\subset\mathbb{C}^4$ be defined by
$$x^2+y^2z+z^{2k-1}+t^{2k-1}=0\,,\quad k\geq 2\,.$$
We see that $X$ is of type $cD_{2k}$. Let us make the blowup 
$\sigma=(k,k-1,1,1)$ (Theorem \ref{T:description}, (i) 1)). It is
easy to see that the exceptional divisor $E$ is given in
$\mathbb{P}(k,k-1,1,1)$ by
$$y^2z+z^{2k-1}+t^{2k-1}=0$$
and that the discrepancy $a(E,X)=1$. The surface $E$ is a cone over
the smooth hyperelliptic curve $y^2z+z^{2k-1}+t^{2k-1}=0$$\subset$
$\mathbb{P}(k-1,1,1)$ of genus $k-1$. This shows that genus
$g(C)$ (see Theorem \ref{T:description}) can be arbitrary.
\end{example}
\begin{example}
Let $X\subset\mathbb{C}^4$ be given by
$$x^2+y^3+yz^3+t^9=0\,.$$
This is a $cE_7$-singularity. Consider the weighted blowup
$\sigma=(5,3,2,1)$ (Theorem \ref{T:description} (iii) 3)). Its
exceptional divisor $E=$$\{y^3+yz^3+y^9=0\}$$\subset$
$\mathbb{P}(5,3,2,1)$ is a cone over a smooth curve of genus $3$.
It is not difficult to show that this curve is non-hyperelliptic.
Indeed, first note that it is trigonal. If it were hyperelliptic,
then we could embed it into $\mathbb{P}^1\times\mathbb{P}^1$ as a
divisor of bidegree $(3,2)$. But then its genus would be
$g\leq 2$.
\end{example}
\begin{example}
Let $X\subset\mathbb{C}^4$ be defined by
$$x^2+y^3+z^5+t^{15}=0\,.$$
It is clear that $X$ is of type $cE_8$. Blowup it with weights 
$(8,5,3,1)$ (Theorem \ref{T:description}, (iv) 6)). Then we obtain
the exceptional divisor 
$$\{y^3+z^5+t^{15}=0\}\subset\mathbb{P}(8,5,3,1)\,.$$
It is a cone over a non-hyperelliptic curve of genus $4$.
\end{example}
\begin{example}
Consider the singularity $X\subset\mathbb{C}^4/\mathbb{Z}_2(0,1,1,1)$
given by
$$x^2+y^3+z^4y+t^4y+z^2t^2+z^6+t^6=0\,.$$
This is a singularity of type $cE/2$. Make the blowups $\sigma_1=
\frac{1}{2}(3,2,3,1)$ and $\sigma_2=\frac{1}{2}(3,2,1,3)$ (see
\cite{Hay}, \S 10). The corresponding exceptional divisors are
\begin{align*}
&E_1\colon\quad\{x^2+y^3+t^4y+t^6=0\}\subset\mathbb{P}(3,2,3,1)
\text{ and }\\
&E_2\colon\quad\{x^2+y^3+z^4y+z^6=0\}\subset\mathbb{P}(3,2,1,3)\,.
\end{align*}
They are cones over elliptic curves. The discrepancies $a_1=a(E_1,X)=
a_2=a(E_2,X)=1/2$. This shows that hyperquotient terminal 
singularities can have more than $1$ nonrational exceptional divisor
with discrepancy $\leq 1$.
\end{example}

\end{document}